\documentclass[11pt, epsfig]{article}
\usepackage{epsfig, amsmath, amssymb, amsthm, times}
\newtheorem {thm}{Theorem}[section]
\newtheorem {lem}[thm]{Lemma}

\theoremstyle{defintion}
\newtheorem {df}[thm]{Definition}

\theoremstyle{remark}
\newtheorem{rem}[thm]{Remark}

\theoremstyle{example}
\newtheorem{ex}[thm]{Example}

\def\qed{\hfill $\Box$ \hfill \\}

\def\R{{\mathbb R}}
\def\Na{{\mathbb N}}

\def\lbl{\label}
\def\be{\begin{equation}}
\def\ee{\end{equation}}
\def\lbl{\label}

\def\pf{\noindent{\em Proof:\ }}
\def\var{{\rm var}}

\def\E{{\mathbb E}}
\def\P{{\mathbb P}}

\def\eps{\epsilon}
\def\veps{\varepsilon}

\def\gs{\sigma}

\def\gl{\lambda}
\def\gk{\kappa}

\title{Stability of equilibria of randomly perturbed maps}

\author{
Pawe{\l } Hitczenko and
Georgi S. Medvedev\thanks{
Department of Mathematics, Drexel University, 3141 Chestnut Street,
Philadelphia, PA 19104,
{\tt phitczen@math.drexel.edu, medvedev@drexel.edu} }
}

\begin{document}
\maketitle
\begin{abstract}
We derive a sufficient condition for stability in probability of an equilibrium 
of a randomly perturbed map in $\R^d$. This condition can be used to stabilize
 weakly unstable equilibria by random forcing. Analytical results on stabilization 
 are illustrated 
with numerical examples of randomly perturbed linear and nonlinear
maps in 
one- and two-dimensional 
spaces.
\end{abstract}

\section{Introduction}
\lbl{sec.intro}
The idea of stabilizing unstable equilibria of dynamical systems by noise 
originates from the pioneering work of Khasminskii on stochastic stability
in the nineteen-sixties \cite{Kha-StochStability}. Stochastic stabilization
has important implications for control theory \cite{Arn90, Mao94, AppMao05, AMR06}
and for numerical methods for stochastic differential equations
\cite{SaiMit96, SaiMit02,Hig00,HMY07,BucKel10}. Furthermore, the 
interplay of stability and noise is important for understanding many dynamical phenomena
in applied science including stochastic synchronization 
\cite{AVR86,DRR06,MP08, GP05},  stochastic resonance  
\cite{Long10, LaiLor10, Fre01},  and noise-induced dynamics
\cite{BerGen06, DVM05, HM09}.

To illustrate the mechanism of stabilization in discrete setting, we
consider a scalar difference equation
\be\lbl{1d}
x_n=(1+\epsilon+\sigma\xi_n) x_{n-1},\; n\in \Na,
\ee
where $0<\epsilon,\sigma\ll 1$ and $(\xi_n)$ are independent copies of the random
variable (RV) $\xi$ with zero mean and $\E\xi^2 =1$. Further,
assume $\P(|\xi| > M)=0$ for some $M>0$. The last condition is used
to simplify the analysis. It can be replaced by a much weaker
condition. For instance, it suffices to have $\E|\xi|^3<\infty$.

For a given $x_0\in\R$, we have
$$
x_n=\left(\prod_{k=1}^n (1+\eps+\sigma\xi_k)\right) x_0.
$$
Let $0<\sigma<M^{-1}$. Then with probability $1$
$$
1+\eps+\sigma\xi_k>0 \; \forall k\in\Na
$$
and we have
$$
\log |x_n| =\log |x_0| +\sum_{k=1}^n \log (1+\eps +\sigma\xi_k)
$$
holding almost surely. By the Strong Law of Large Numbers,
$$
n^{-1}\sum_{k=1}^n \log (1+\eps +\sigma\xi_k)\to \E\log (1+\eps +\sigma\xi) \;\mbox{as}\; n\to\infty
$$
almost surely.
Thus, the asymptotic stability of the origin (in the almost sure sense)  will follow if
\be\lbl{neg-log}
\E\log (1+\eps +\sigma\xi)<0.
\ee
Using the Taylor expansion of $\log$ and $\E\xi=0$, we have
$$
\E\log (1+\eps+\sigma\xi)=\eps-{\sigma^2\over 2}+O(\sigma^3,\eps^3).
$$
Thus, the stabilization of the weakly unstable equilibrium of \eqref{1d} is achieved if
\be\lbl{stab-1d}
\eps-{\sigma^2\over 2} <0
\ee
for $0<\eps,\sigma\ll 1$. A similar stabilization condition is known for ordinary 
differential equations \cite{Kha-StochStability}.

The stability analysis for \eqref{1d} can be extended to linear maps in $\R^d$ using 
the Furstenberg-Kesten theory \cite{FurKes60}. 
For scalar nonlinear difference equations, stabilization was studied by Appleby, Mao, and Rodkina 
\cite{AMR06} and by Appleby, Berkolaiko, and Rodkina \cite{ABR09} (see
also \cite{ABR08, BR12, AKMR10, BR06}).
Certain higher-dimensional models were analyzed  in the
context of stability of finite-difference schemes (see \cite{BucKel10} and references therein).
In this paper, we show that one can achieve stabilization
with high probability  for a general $d$-dimensional nonlinear map 
 under fairly general assumptions on the stabilizing perturbation.
Specifically, we study the following difference equation in $\R^d$
\be\lbl{model}
x_{n+1}=(A+B)x_n +q(x_n),
\ee
where $q(x)=O(|x|^2)$ is a smooth function, $A$ and $B$ are deterministic and stochastic
$d\times d$ matrices respectively.  We assume that the spectral radius of $A$ is slightly
greater than $1$, $\rho(A)=1+\eps,$ $0<\eps\ll 1$ and ask how to choose mean-zero matrix 
$B=B(\eps)$ to stabilize the equilibrium at the origin. Our motivation for 
considering (\ref{model}) is two-fold. On one hand, we want to
understand how to tame weak instability in general $d$-dimensional
maps by noise. Eventually, we want to apply these results
to stabilize periodic orbits of randomly perturbed stochastic ordinary
differential equations in
$\R^{d+1}$. In this case, (\ref{model})
represents a Poincare map \cite{HM13}.  Stochastic stabilization 
of period orbits remains largely unexplored area of research
with many promising applications.

The organization of this paper is as follows. In the next section, we prove a sufficient condition
for stability (in probability) of an equilibrium in a $d$-dimensional map 
(cf. Theorem~\ref{thm.Lyap}). To prove this theorem, we use the 
Strong Law of Large Numbers to show that the Lyapunov exponent of a 
typical trajectory is negative. The rest of the proof follows an
argument developed 
for deterministic dynamical systems \cite{KocPal10}. In \S \ref{sec.stabilization}
we apply Theorem~\ref{thm.Lyap} to the problem of stabilization. In \S\ref{sec.numerics},
we illustrate our results with several numerical examples using one- and two-dimensional systems.

\section{Stochastic stability} \lbl{sec.stability}

Consider an initial value problem for the following difference equation
\be\lbl{nmap}
x_n=M_nx_{n-1}+q(x_{n-1}), \; n\ge 1.
\ee
where $(M_n)$ are independent copies of a $d\times d$ random matrix $M$;
$q:\R^d\to\R^d$ is a continuous function such that
\be\lbl{quad}
|q(x)|\le C_1|x|^2, \; x\in B_\delta=\{x: |x|\le\delta\}
\ee
for some $C_1,\delta>0$. Here and below, we will use $|\cdot|$ to
denote the Euclidean norm of a vector.
The initial condition $x_0$ is assumed to be deterministic.

\begin{df} (cf.~\cite{Kha-StochStability}) The equilibrium at the origin of (\ref{nmap})
is said to be stable in probability if for any $\veps>0$
$$
\lim_{|x_0|\to 0} \P\{ \sup_{n\ge 1} |x_n|>\veps\}=0.
$$
\end{df}

\begin{thm}\lbl{thm.Lyap}
Suppose
\be\lbl{Lyapunov}
0<\lambda=-\E\log\|M\|<\infty.
\ee
Then the equilibrium at the origin of (\ref{nmap}) is
stable in probability.
\end{thm}
\begin{rem}\lbl{rem.norm}
In (\ref{Lyapunov}), $\|\cdot\|$ is an arbitrary matrix norm. 
The same matrix norm is used throughout this section.
\end{rem}

Condition (\ref{Lyapunov}) guarantees that the largest Lyapunov exponent of 
a generic trajectory is negative. This implies stability of $x_n\equiv 0$ with high
probability. Theorem~\ref{thm.Lyap} is a stochastic counterpart of the 
result of Ko{\c{c}}ak and Palmer for deterministic maps
\cite[Theorem~4]{KocPal10}.
It follows immediately from the proof of the following lemma,
which also yields the rate of convergence of $(x_n)$ to the origin.

\begin{lem}\lbl{thm.stability}
Let $(x_n)$ denote a trajectory of \eqref{nmap} subject to \eqref{Lyapunov}.
Then for any $0<\veps<\min\{1,\lambda/3\}$ there exist $\eta>0,$
$\delta_1>0,$ and $\mu=\exp\{-\lambda+\veps\}<1$ such that
\be\lbl{stability}
|x_i|\le \eta\mu^i,\; i=0,1,2,\dots
\ee
with probability at least $1-\veps$ provided $|x_0|\le \delta_1.$ 
\end{lem}
\pf
Suppose $0<\veps<\min\{1,\lambda/3\}$ is arbitrary but fixed.
Let $\lambda_k:=\log \|M_k\|$ and note that 
\[
\frac1n\sum_{k=1}^n\lambda_k\stackrel{a.s.}\longrightarrow -\lambda<0\;\mbox{as}\; n\to\infty,
\]
by the Strong Law of Large Numbers \cite[Theorem~22.1]{bil}.
Thus, there exists $n_0>1$ such that 
$$
\P\left( \bigcup_{n\ge  n_0}\left\{ 
\left|\frac1n\sum_{k=1}^{n}\lambda_{k}+\lambda \right|
>\veps
\right\}
\right)
<\frac\veps2,
$$
i.e., for $n\ge n_0,$
\be\lbl{exp}
-\lambda-\veps\le \frac1n\sum_{k=1}^{n}\lambda_{k}\le -\lambda+\veps
\ee
holds on the set of probability at least $1-\veps/2$. 
In the remainder of the proof, we restrict to the realizations $(M_k)$ for which (\ref{exp}) holds.

Using (\ref{exp}), for any $n\ge k\ge n_0>1$, we have
\be\lbl{product-1} 
\prod_{j=k}^n\|M_j\|=\frac{\prod_{j=1}^n\|M_j\|}{\prod_{j=1}^{k-1}\|M_j\|}\le
\exp\{n(-\lambda+\veps)-(k-1)(-\lambda-\veps)\}=
\mu^{n-k+1}e^{2(k-1)\veps}.
\ee
Similarly, for every $ 1\le k<n_0,$ we have
\begin{align}\nonumber
\prod_{j=k}^n\|M_j\|&=\left(\prod_{j=k}^{n_0-1}\|M_j\|\right)\left(\prod_{j=n_0}^n\|M_j\|\right)\le
\left(\prod_{j=k}^{n_0-1}\|M_j\|\right)\mu^{n-n_0+1}e^{2(n_0-1)\veps}\\
\lbl{product-2}
&\le \bar M_{n_0} \mu^{n-k+1}e^{2(k-1)\veps},
\end{align}
where 
$$
\bar M_{n_0}= \max_{1\le k\le n_0-1}\left\{ \mu^{k-n_0} e^{2(n_0-k)\veps} \prod_{j=k}^{n_0-1}\|M_j\|  \right\}.
$$

Since $\bar M_{n_0}$ is an
integrable random variable, by Markov inequality, we have 
$$
\P\left( \bar M_{n_0}\ge M\right)\le {\E\bar M_{n_0}\over M} \quad\forall M>0.
$$
Choosing $M=M(\veps)>0$ sufficiently large, we have 
\be\lbl{Markov}
\P\left( \bar M_{n_0} \ge M \right)\le\frac\veps2.
\ee
The combination of \eqref{product-1}, \eqref{product-2}, and \eqref{Markov}
yields 
\be\lbl{prod}
 \prod_{j=k}^n\|M_j\|\le C_2\mu^{n-k+1} e^{2(k-1)\veps}, \quad 1\le k\le n,
\ee
holding with probability at least $1-\veps$, where $C_2=\max\{M,1\}$ depends 
on $\veps$ but not on $n$ or $k$.

We are now in a position to prove \eqref{stability}. To this end, fix $0<\eta\le \delta$ and
choose $0<\delta_1\le \eta$ such that
\be\lbl{choose-delta1}
C_2\delta_1 \exp\{ C_1C_2\eta/ 1-\nu\} \le \eta,
\ee
where $\nu:=e^{-\lambda+3\veps}<1$. With these constants $\eta$ and $\delta$,
we will show \eqref{stability} by induction.

The claim in \eqref{stability} obviously holds for $i=0.$ Let $p\ge 1$ and suppose that
\be\lbl{hypothesis}
|x_i|\le \eta \mu^i
\ee
holds for $i=0,1,\dots, p-1.$ We want to show that this entails 
$$
|x_p|\le \eta \mu^p.
$$
Iterating (\ref{nmap}), we have
\be\lbl{pstep}
x_p=\left(\prod_{k=0}^{p-1} M_{p-k}\right) x_0 + \sum_{j=1}^{p}
\left(\prod_{k=0}^{p-j} M_{p-k}\right) q(x_{j-1}).
\ee
Using the triangle inequality, submultiplicativity of the matrix norm, and
 (\ref{quad}), from (\ref{pstep}) we obtain
$$
|x_p|\le \left(\prod_{k=0}^{p-1} \|M_{p-k}\|\right) |x_0| + C_1 \sum_{j=1}^{p}
\left(\prod_{k=0}^{p-j} \|M_{p-k}\|\right) |x_{j-1}|^2.
$$
Here, we also used the induction hypothesis \eqref{hypothesis}, which implies that 
$|x_j|\le \delta,\, j=0,1,\dots, p-1$ so that (\ref{quad}) is applicable.
Using  (\ref{prod}), we further derive
$$
|x_p|\le C_2\mu^p |x_0| + C_1C_2\sum_{j=1}^{p} \mu^{p-j+1} e^{2(j-1)\veps} |x_{j-1}|^2.
$$
Using the induction hypothesis \eqref{hypothesis},  we continue
\be\lbl{done-xp}
|x_p|\le C_2\mu^p |x_0| + C_1C_2\eta \mu^{p}\sum_{j=1}^{p} e^{2(j-1)\veps} |x_{j-1}|.
\ee
Next, we rewrite \eqref{done-xp} in terms of
\be\lbl{def-z}
z_i=\mu^{-i} |x_i|, \; i=0,1,2,\dots,p,
\ee
to obtain
$$
z_p\le C_2 z_0 +C_1C_2\eta\sum_{j=1}^{p} \nu^{j-1}z_{j-1}, \quad \nu=e^{-\lambda+3\veps}<1.
$$
By the discrete Gronwall's inequality (see ~Lemma~\ref{lem.Gronwall} below),
we have
$$
z_p\le C_2 z_0 \exp\left\{C_1 C_2\eta\sum_{k=1}^p \nu^{k-1}\right\}\le C_2 \delta_1 
\exp\{C_1C_2\eta/(1-\nu)\}\le \eta,
$$
where we used \eqref{choose-delta1} in the last inequality.
Recalling the definition of $z_p$ \eqref{def-z}, we conclude that $|x_p| \le \eta\mu^p.$ 
\\ \qed
 
\begin{lem}(cf.~\cite{KocPal10})
\lbl{lem.Gronwall}
Let $\{z_k\}_{k=0}^\infty$and $\{\mu_k\}_{k=1}^\infty$ be two nonnegative
sequences such that
\be\lbl{growth-Gron}
z_k\le B+\sum_{j=1}^{k}\mu_j z_{j-1}, \; k\in [p],
\ee
for some $p\in\Na$. Then for $k\in [p]$
$$
z_k\le B\exp\left\{\sum_{j=1}^k \mu_j\right\}.
$$
\end{lem}

\section{Stabilization}\lbl{sec.stabilization}
\setcounter{equation}{0}
Consider the following difference equation in $\R^d$:
\be\lbl{map}
x_{n+1}=(A(\eps)+B_n(\eps))x_n+f(x_n), \; n=0,1,2,\dots,
\ee
where $f(x)=O(|x|^2)$ and $A(\eps)\in \R^{d\times d}$ is a matrix with the 
spectral radius
\be\lbl{radius}
\rho(A(\eps)):=\max\{|\gl|:\
\gl\  \mbox{is an eigenvalue of $A(\eps)$}\}=1+\eps, \;\; 0<\epsilon\ll 1.
\ee
$(B_n(\eps))$ are independent copies of a random matrix
$B(\eps)\in\R^{d\times d}$.
We want to identify the conditions on $B(\eps)$, which guarantee
stabilization of the unstable equilibrium at the origin. To keep the notation simple, 
we will freely suppress the dependence of $A$ and $B$ on $\eps$, whenever it is 
not essential.

Suppose the  Jordan normal form of $A$ is 
\be\lbl{Jordan}
A^0+U,
\ee
where $A^0$ is the block-diagonal matrix
\be\lbl{block-diag}
A^0=\operatorname{diag}\left( A_1, A_2, \dots, A_k\right)
\ee
for some $k\in [d]$. Block $A_i, \; i\in [k],$ is $(\lambda_i)$ if the 
corresponding eigenvalue of $A$ is real, or 
$$
\begin{pmatrix} a_i & -b_i \\b_i & a_i\end{pmatrix},
$$
otherwise.
Therefore, 
\be\lbl{norm-A0}
\|A^0\|=\rho(A)=1+\epsilon.
\ee
Here and below, $\|\cdot\|$ stands for the operator norm of a matrix.

The upper-triangular matrix $U$ is nonzero only if $A$ has multiple eigenvalues.
In this case, it has the following form
$$
U=\begin{pmatrix} 
O_{11} &*&*& \dots&*\\ 
 & \dots & *&\dots& *\\
 & \dots & &\dots&  &\\
0 &   & \dots &  & O_{kk}
\end{pmatrix},
$$
where $O_{ii}$ is a $d_i\times d_i$ zero block whose dimension coincides with
that  of $A_i$ for each $i\in [k]$. By changing coordinates, one can achieve
\be\lbl{small-U}
\|U\|<\kappa
\ee
for any $\kappa>0$ given in advance.\footnote{Indeed, let
  $D_t=\operatorname{diag}\left(tI_1,t^2 I_2, \dots, t^k I_k\right), t>0,$ where $I_i, \, i\in [k],$
is a $d_i\times d_i$ identity matrix.
Then all entries above the main diagonal of $D_t UD_t^{-1}$ can be made arbitrarily small provided 
$t$ is large enough.}

Thus, without loss of generality, we assume that matrix $A$ 
in \eqref{map} has the following
form
\be\lbl{normal-A}
A=A^0+U,
\ee
where the block-diagonal matrix $A^0$ and the upper diagonal matrix $U$ are subject to \eqref{block-diag},
\eqref{norm-A0}, and \eqref{small-U}, respectively.

Next, we formulate our assumptions on the random matrix $B$. First, we
describe a general class of stabilizing random matrices. Later, we
will see that in practice stabilization can be achieved with a very simple
random matrix $B(\eps)$.

Let 
\be\lbl{take-B}
B(\eps)=A(\eps)G(\eps), 
\ee
where $G(\eps)$ is a $d\times d$ symmetric matrix, whose entries $g_{ij}(\eps)$ are mean zero
non--degenerate RVs  with finite third moments subject to the following conditions. 

Denote
$$
\gs_{ij}^2(\eps):= \E g_{ij}(\eps)^2, \; (i,j)\in [n]^2, \quad \mbox{and}\quad  
\sigma(\eps):=(\sigma_{11}(\eps),\sigma_{22}(\eps),\dots, \sigma_{nn}(\eps)).
$$
We assume
\begin{eqnarray}
\lbl{vanishing-sigma}
&& \lim_{\eps\to 0} |\sigma(\epsilon)|=0, \\
\lbl{diagonal}
&&  \lim_{\eps\to 0}{\gs_{ij}(\eps)\over \gs_{ii}(\eps)^2}=0, \; i\neq j,\\
\lbl{3rd_mom}
&&  (\E |g_{ij}(\eps)|^3)^{1/3}\le K \gs_{ij}(\eps),\quad 1\le i,j\le d,
\end{eqnarray}
for some $K>0$ independent of $\eps$.
\begin{rem}
Condition \eqref{3rd_mom}  is  easy to fulfill. For example,
we may take $g_{ij}(\eps)=\eps\xi_{ij}$, where $\xi_{ij}$ is a mean--zero
random variable with the finite third moment, $1\le i,j\le d$. Then, for
each such $i$ and $j$, $(\E|\xi_{ij}|^3)^{1/3}\le
K_{ij}(\E\xi_{ij}^2)^{1/2}$ for some constant $K_{ij}$ and \eqref{3rd_mom}
holds with $K=\max_{1\le i,j\le d}K_{ij}$  for all $\eps>0$.
In particular, if $(\xi_{ij})$ are (arbitrarily dependent) standard normal
random variables then \eqref{3rd_mom} holds with $K=2\sqrt{2/\pi}$.
\end{rem}

By Theorem~\ref{thm.Lyap}, for stabilization in \eqref{map} it is sufficient to show that the condition 
\be\lbl{Kesten}
\E\log\|A(\eps)+B(\eps)\| <0
\ee
holds for some small $\eps>0$. The following lemma provides a sufficient condition for \eqref{Kesten}.

\begin{lem}\lbl{lem.diagonal} 
Suppose $A(\eps)$ and $B(\eps)$ satisfy the assumptions 
\eqref{norm-A0}-\eqref{3rd_mom}.
Then (\ref{Kesten}) holds for sufficiently small $\epsilon>0$,
provided
\be\lbl{halfsq}
1<\limsup_{\eps\searrow 0} {|\sigma(\epsilon)|^2\over 2\eps}<\infty.
\ee
\end{lem}
\begin{rem}\lbl{rem.parametric}
The parametric dependence $B(\eps)$ in \eqref{map} is used for convenience of presentation
only. By interpreting $\sigma$ as a function of $\eps$, we are dealing with a single small parameter
$\epsilon$, instead of having to work with both $\eps$ and $|\sigma|$. The parametric dependence
in \eqref{halfsq} is not essential. What this condition means is that $|\sigma|^2$ should be 
large enough compared to $\eps$, while both $|\sigma|$ and $\eps$ must be small.
\end{rem}

\pf By the submultiplicativity of the matrix norm and \eqref{take-B}, we have
\be\lbl{logA+B}
\log\|A+ B\|=\log\|A(I+ G)\|\le\log\|A\|+\log\|I+ G\|.
\ee
Let $\eps$ and $\gk=\gk(\eps)$ whose values will be specified later be chosen. Using \eqref{normal-A} and \eqref{norm-A0}, from \eqref{logA+B} we further
obtain
\be\lbl{bd4log}
\log\|A+ B\| \le\log (1+\eps+\gk)+\log \|I+ G\|.
\ee

By Gershgorin Theorem (cf.~\cite{Horn-Johnson}),
\[
\|I+ G\|=
\rho(I+ G)  \le\max_i\Big(|1+ g_{ii}|+ \sum_{j\ne i}|g_{ij}|\Big).
\]
By the monotonicity of logarithm,
\[
\log\|I+ G\|\le\max_i\log (|1+ g_{ii}|+\sum_{j\ne i}| g_{ij}|)\le \sum_i 
\log(|1+ g_{ii}|+\sum_{j\ne i}|g_{ij}|).
\]   
Taking expectations on both sides, we get 
\[
\E\log \|I+ G\| \le \sum_{i=1}^d \E\log \Big( |1+ g_{ii}|+\sum_{j\ne i}|g_{ij}| \Big).\]
For each $i$ 
\begin{eqnarray} \nonumber
\E\log\Big(|1+
  g_{ii}|+\sum_{j\ne i}|g_{ij}| \Big)&\le& 
\E\log(1+ g_{ii} +\sum_{j\ne i}|g_{ij}|)I_{|g_{ii}|<1}\\&+&
\E\log\Big(1+| g_{ii}|+\sum_{j\ne i}|g_{ij}|\Big)I_{|g_{ii}|\ge1}.\label{2nd_log}
\end{eqnarray}
By expanding the logarithm in the first term and using the fact that
$\E g_{ii}=0$ we get 
\begin{eqnarray}
&&\nonumber
 \E\left(g_{ii}+\sum_{j\ne i} |g_{ij}| -\frac{(g_{ii}+\sum_{j\ne i}|g_{ij}| )^2}2+O(
(|g_{ii}|+\sum_{j\ne i}|g_{ij}| )^3)\right)I_{|g_{ii}|<1}\\&&\nonumber\quad=
\E\left(g_{ii}+\sum_{j\ne i} |g_{ij}| -\frac{1}2(g_{ii}+\sum_{j\ne i}|g_{ij}| )^2+O(
(|g_{ii}|+\sum_{j\ne i}|g_{ij}| )^3)\right)\\&&\nonumber\qquad -
\E\left(g_{ii}+\sum_{j\ne i} |g_{ij}| -\frac{1}2(g_{ii}+\sum_{j\ne i}|g_{ij}| )^2+O((
|g_{ii}|+\sum_{j\ne i}|g_{ij}| )^3)\right)I_{|g_{ii}|\ge1}\\
&&\nonumber\quad=
\sum_{j\ne i}\E|g_{ij}| -\frac12\E g_{ii}^2-\sum_{j\ne i}\E g_{ii}|g_{ij}|-\frac12\E(\sum_{j\ne i}|g_{ij}| )^2+
\E O((|g_{ii}|+\sum_{j\ne i}|g_{ij}| )^3)\\
\lbl{3.xx}
&&\qquad
 +O(\sum_{m=1}^3\E(|g_{ii}|+\sum_{j\ne i}|g_{ij}| )^mI_{|g_{ii}|\ge1}).
\end{eqnarray}
Note that since $\log(1+x)\le x$ for $x\ge 0$, the bound on the last term on the right--hand side
of \eqref{3.xx} gives the bound for the second term on 
the right--hand side in  \eqref{2nd_log}. 

We estimate the
terms above as follows
\begin{eqnarray*}\sum_{j\ne i}\E|g_{ij}|&=&\sum_{j\ne i}O(\gs_{ij})=o(\gs^2_{ii}),
\quad(\mbox{by \eqref{diagonal}})\\
\big|\sum_{j\ne i}\E g_{ii}|g_{ij}|\big|&\le&\gs_{ii}\sum_{j\ne i}
\gs_{ij}=o(\gs^2_{ii}),\quad(\mbox{by the Cauchy-Schwarz inequality and \eqref{diagonal}})\\
\E (\sum_{j\ne i}|g_{ij}| )^2&=&\sum_{j\ne i}O(\gs_{ij}^2)=o(\gs^2_{ii}),\quad(\mbox{by \eqref{diagonal}})\\
\E (|g_{ii}|+\sum_{j\ne i}|g_{ij}| )^3&=&O(\E |g_{ii}|^3)+\sum_{j\ne i}O(\E |g_{ij}|^3)\\  
\E (|g_{ii}|+\sum_{j\ne i}|g_{ij}| )^m I_{|g_{ii}|>1}&=& O(\E |g_{ii}|^mI_{|g_{ii}|>1})+\sum_{j\ne i}O(\E |g_{ij}|^m).
\end{eqnarray*}
For $m=1,2$ and $j\ne i$, $\E |g_{ij}|^m=o(\gs_{ij}^2)$ as verified above. Further,  for $1\le m\le 3$
\[\E |g_{ii}|^mI_{|g_{ii}|>1}\le \E |g_{ii}|^3I_{|g_{ii}|>1}\le \E |g_{ii}|^3.
\]
Hence, by \eqref{vanishing-sigma}, \eqref{diagonal}, and \eqref{3rd_mom} for all $1\le i,j\le d$, 
\[\E |g_{ij}|^3=O(\gs_{ij}^3)=o(\gs_{ij}^2)=o(\gs_{ii}^2).\]
Plugging all of this into \eqref{bd4log} and using 
 $\log(1+\eps+\gk)\le \eps+\gk$  we obtain that 
\be\lbl{last_bd} \E \log\|A+B\| \le\eps+\gk-\frac12\sum_{i=1}^d(\gs_{ii}^2+o(\gs_{ii}^2))= \eps\left(1-\frac{|\gs(\eps)|^2}{2\eps}+\frac\gk\eps+\frac{o(|\gs(\eps)|^2)}\eps\right)
.\ee

Let $K<\infty$ be any number strictly larger than $\limsup(|\gs(\eps)|^2/(2\eps))$. By  \eqref{halfsq} there exists an $\eps_0>0$ such that for all $0<\eps<\eps_0$ 
\[1+\eps_0\le\frac{|\gs(\eps)|^2}{2\eps}\le K
.\]
Decreasing $\eps$ if necessary we may assume that the error term $o(|\gs(\eps)|^2)$ satisfies
\[\frac{o(|\gs(\eps)|^2)}{|\gs(\eps)|^2}\le\frac{\eps_0}{6K}.
\]
Finally, for the chosen $\eps$ we choose  a $\gk>0$ satisfying $\gk\le \eps\cdot\eps_0/3$. With these choices the right--hand side of \eqref{last_bd} is at most
\[\eps\left(-\eps_0+\frac{\eps_0}3+\frac{o(|\gs(\eps)|^2)}{|\gs(\eps)|^2}\cdot\frac{|\gs(\eps)|^2)}\eps\right)
\le\eps\left(-\eps_0+\frac{\eps_0}3+\frac{\eps_0}{6K}\cdot2K\right)\le-\frac{\eps\cdot\eps_0}3<0.
\]
This proves that the left--hand side of \eqref{last_bd} is negative and completes the proof. 
\qed
\begin{rem}
As can be easily seen from the proof, for stabilization of the unstable equilibrium in \eqref{map}
it is sufficient to take a diagonal matrix $D=\operatorname{diag}(g,g,\dots,g),$ where mean zero RV $g$
meets the conditions on the three first moments \eqref{vanishing-sigma} and \eqref{3rd_mom} as well as \eqref{halfsq}.
In particular, one can take $g=a\xi$, where $\xi$ is a standard normal
RV and $a=a(\eps)\to 0$, but $\lim\limits_{\eps\to 0}(a/\sqrt{2\eps})>1$.
Thus, in practice, it suffices to use a single RV to stabilize a weakly unstable equilibrium in $\R^d$.
 \end{rem}

\begin{figure}
\begin{center}
 {\bf a}\includegraphics[height=1.8in,width=2.0in]{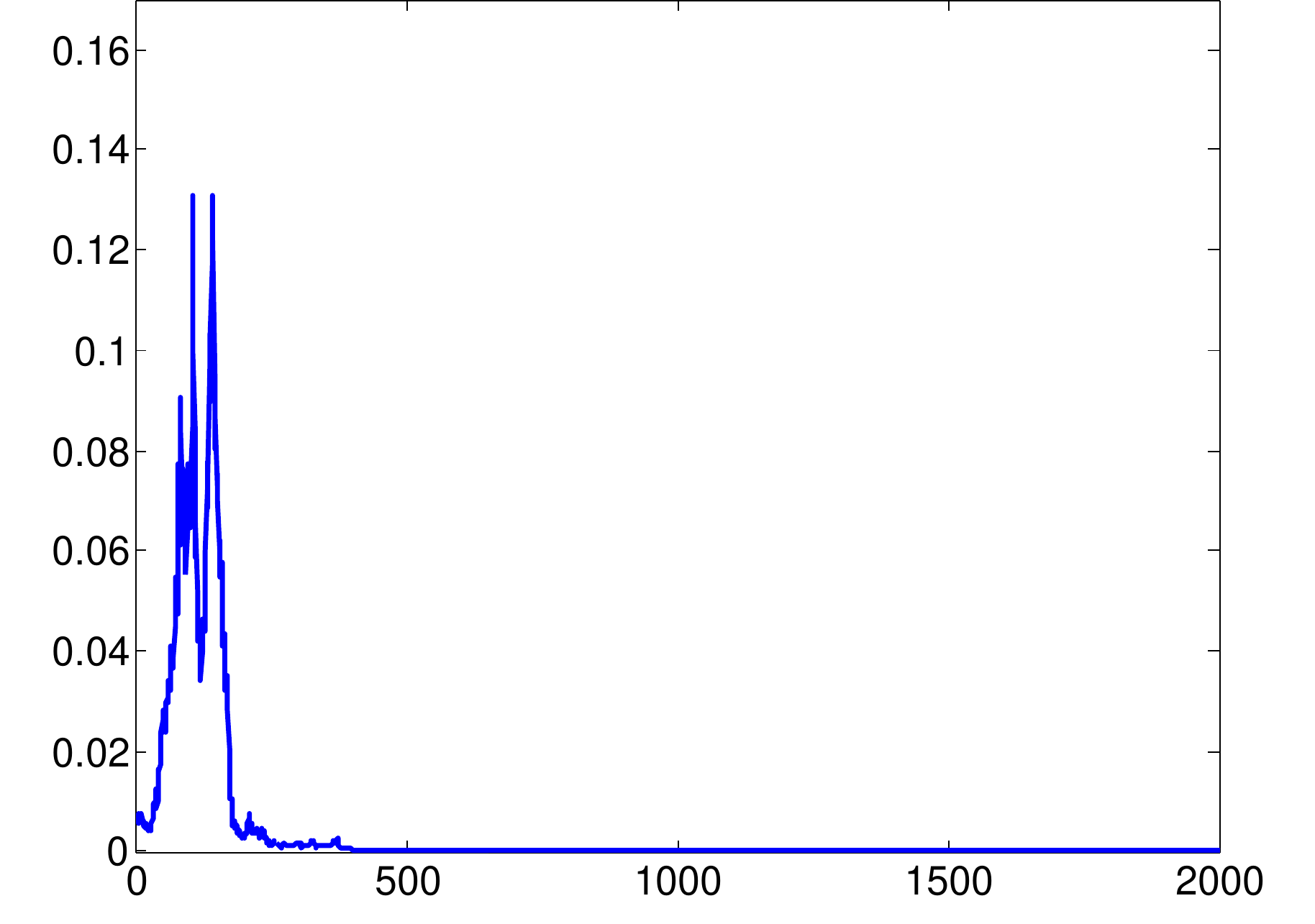}\quad
 {\bf b}\includegraphics[height=1.8in,width=2.0in]{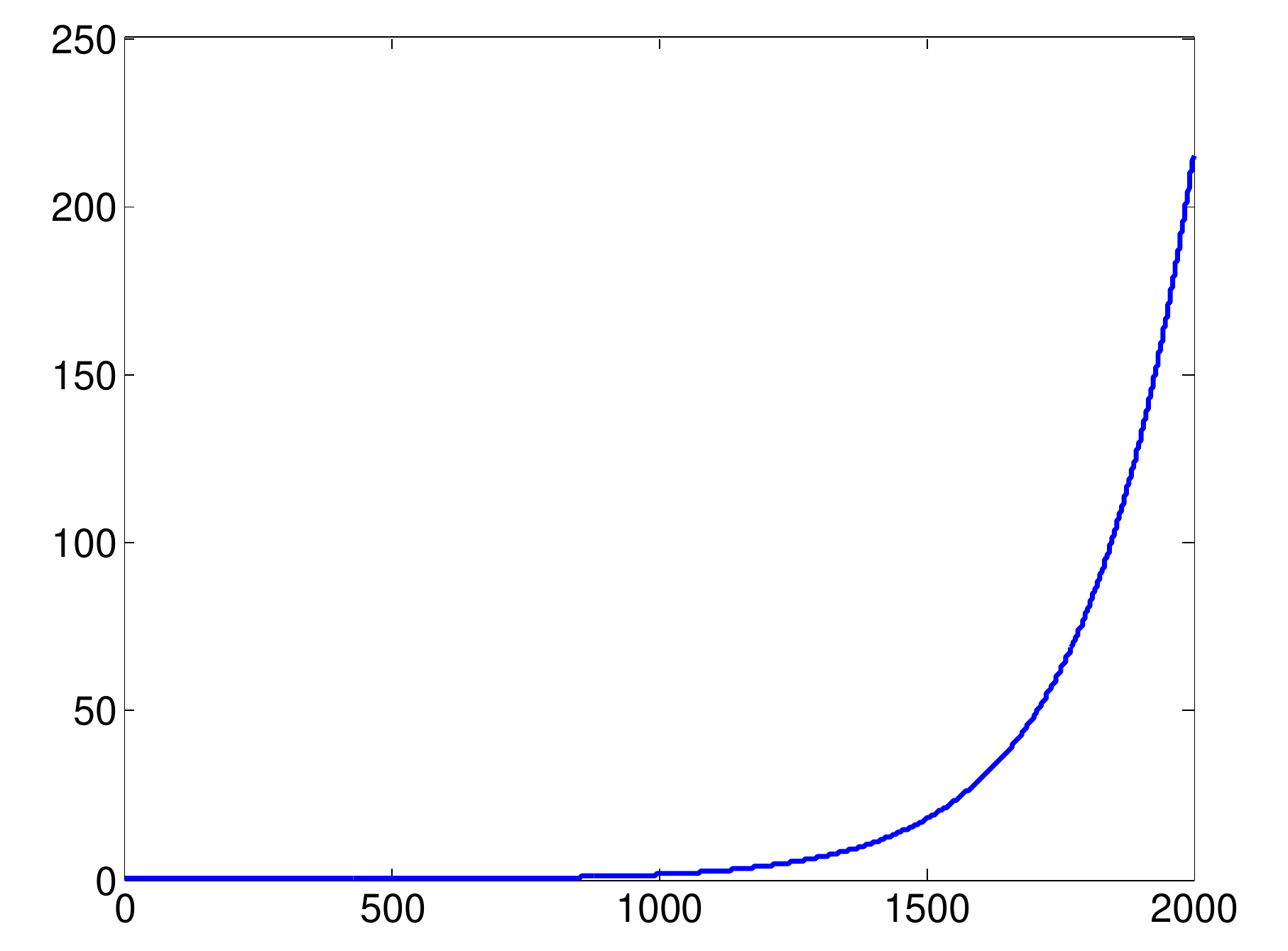}
\end{center}
\caption{ \textbf{a}) Time series generated by the  stochastic  one-dimensional system defined in 
Example~\ref{ex.1d-1}. The values of parameters are $\eps=0.005,$ $\rho=4$.
\textbf{b}) The time series generated by the underlying deterministic system
($\rho=0$) is included for comparison.
}
\lbl{f.1}
\end{figure}

\section{Examples}\lbl{sec.numerics}
\setcounter{equation}{0}
In this section, we illustrate our analysis of stabilization  with
several  numerical examples. 

\subsection{One-dimensional maps}\lbl{sec.1D}
We consider first a scalar difference equation
\be\lbl{scalar}
x_{n+1}=f(x_n)+\xi_{n+1}x_n, \quad n=0,1,2,\dots,
\ee
where $f:\R\to\R$ is a smooth function, $f^\prime (0)=1+\eps$, and
$(\xi_n)$ are independent copies of a RV $\xi$ with $\sigma^2:=\var(\xi)<\infty$.

Lemma~\ref{lem.diagonal} yields 
\be\lbl{stab-scalar}
{\sigma^2\over 2}-\eps>0
\ee
as a sufficient condition for stabilization provided $\eps$ and $\sigma$ are small
enough.

\begin{ex}\lbl{ex.1d-1} Let $f(x)=(1+\eps)x,$ $\sigma^2=\rho\eps,$ and 
$\xi\in\mathcal{N}(0,\sigma^2)$. The results of numerical simulations of (\ref{scalar}) 
with the linear map above
with small positive initial condition
are shown in Figure~\ref{f.1}. Plot \textbf{a} shows that the trajectory of the 
random system with noise intensity subject to (\ref{stab-scalar}) after a brief explosion
converges to the origin. The deterministic trajectory in \textbf{b} grows exponentially. 
\end{ex}

\begin{ex}\lbl{ex.1d-2} Next, we consider a nonlinear map $f(x)=\lambda x(1-x)$. 
For $\lambda=1+\epsilon> 1$,
  the logistic map $f$ has two fixed points: $\bar x_1=0$ and $\bar x_2=\epsilon (1+\epsilon)^{-1}$.
For $0<\epsilon\ll 1$, the former is unstable, while  the latter is stable.
All trajectories of the deterministic map $x\mapsto f(x)$ starting
from $x_0\in (0,1)$ converge to $\bar x_2$ (see Fig.~\ref{f.2}\textbf{b}).
In the presence of noise, however, the
iterations of (\ref{scalar}) with high probability converge to $\bar x_1,$ 
provided (\ref{stab-scalar}) holds
and $\eps$ is small enough (see Fig.~\ref{f.2}\textbf{a}).
\end{ex}

\begin{figure}
\begin{center}
 {\bf a}\includegraphics[height=1.8in,width=2.0in]{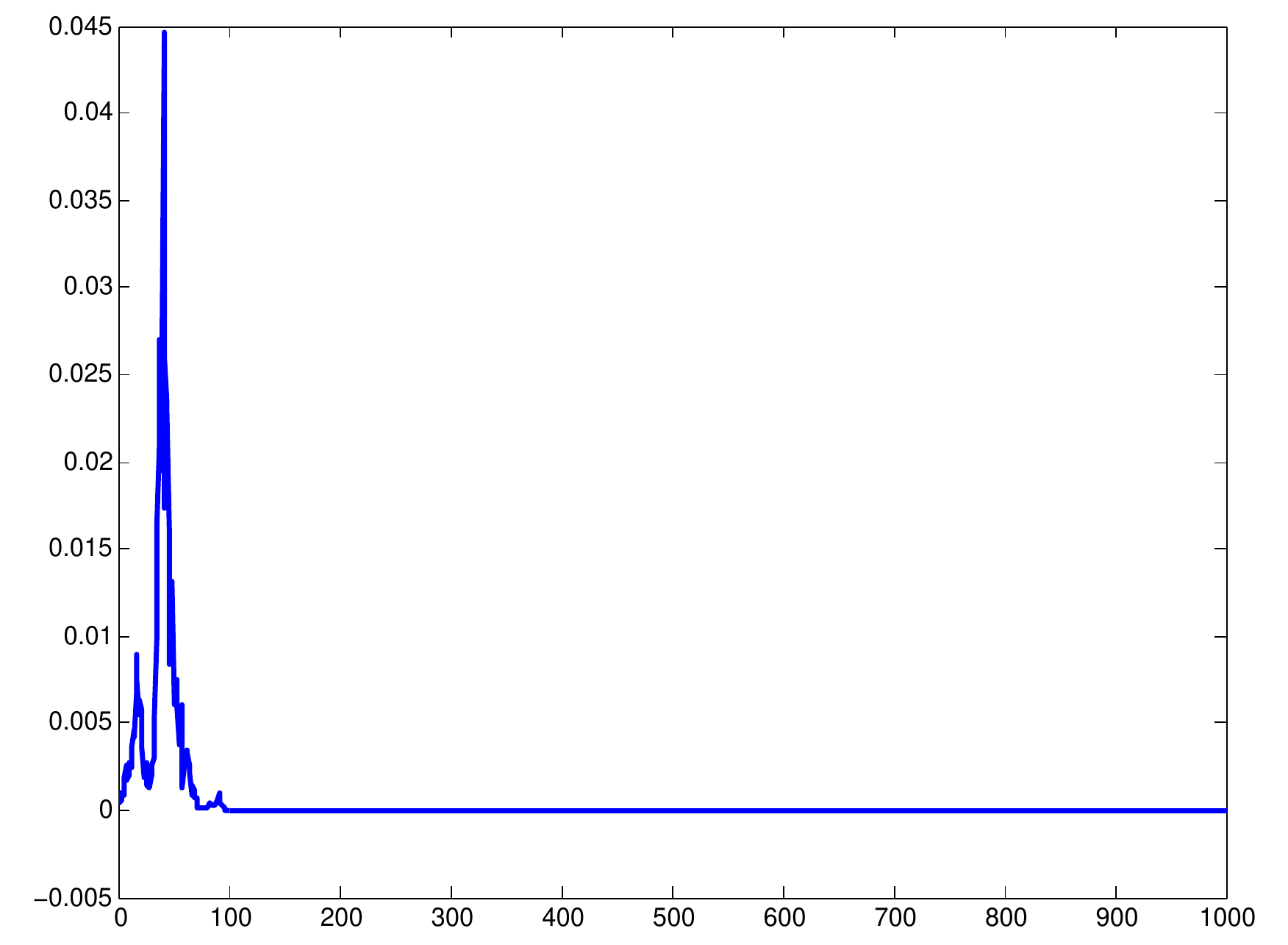}\quad
 {\bf b}\includegraphics[height=1.8in,width=2.0in]{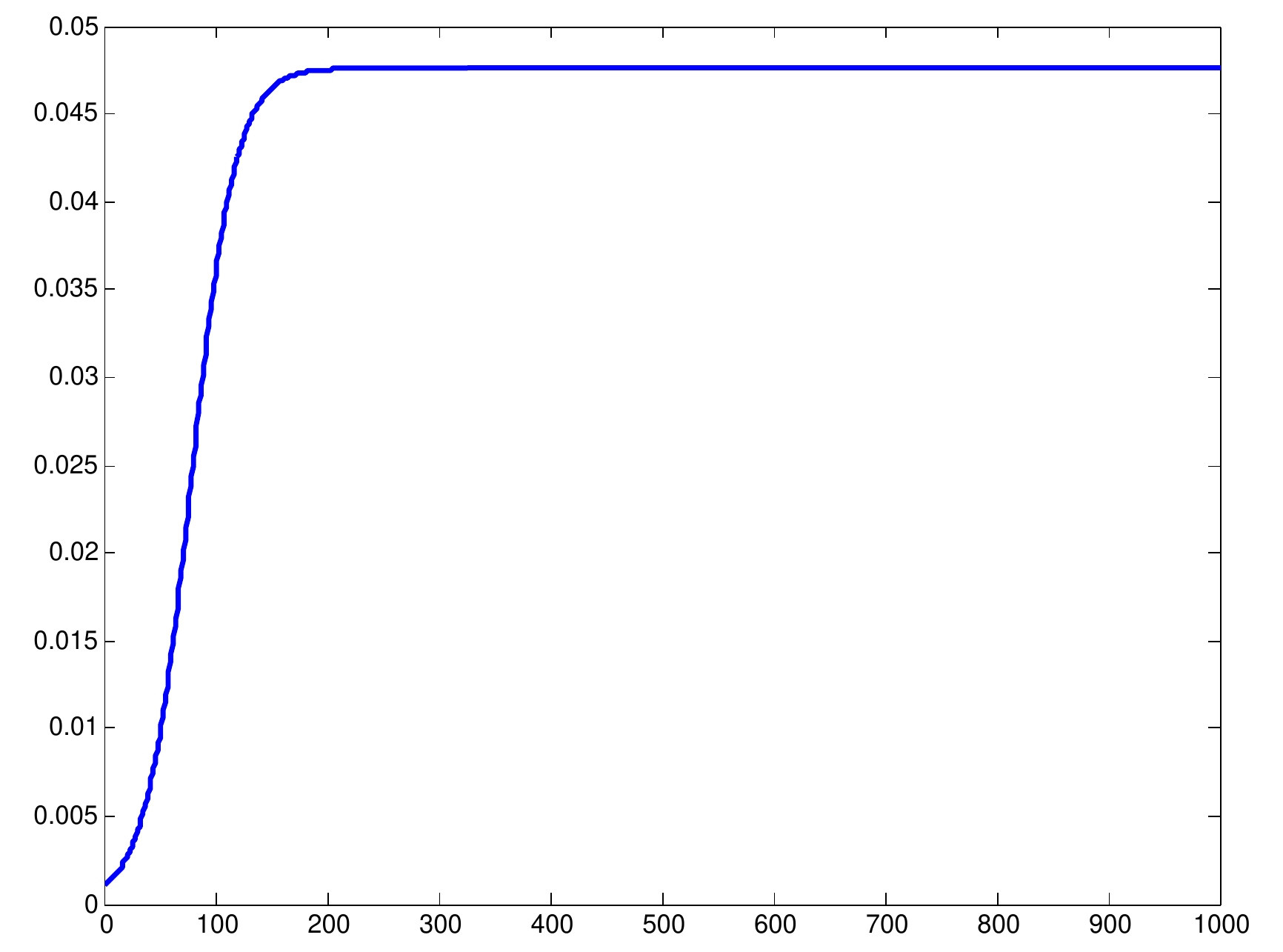}
\end{center}
\caption{ \textbf{a}) Time series generated by the randomly perturbed logistic map  
(see Example~\ref{ex.1d-2}).  Here, $\xi\in\mathcal{N}(0,\rho\eps)$ and
the values of parameters are $\eps=0.05,$ $\rho=3$.
\textbf{b}) The time series generated by the underlying deterministic system
($\rho=0$) is included for comparison.
}
\lbl{f.2}
\end{figure}

\subsection{Two-dimensional maps}\lbl{sec.2D}
We next turn to the $2$D case. To this effect, we consider
\be\lbl{planar}
x_{n+1}=(A+B)x_n, \quad n=0,1,2,\dots,
\ee
where $A$ is a $2\times 2$ deterministic matrix  and
\be\lbl{defB}
B=\sigma \begin{pmatrix}   \xi_{11} & \eps\xi_{12} \\
                                          \eps\xi_{12} & \xi_{22}
                \end{pmatrix},
\quad \xi_{ij}\in \mathcal{N}(0,1), \sigma^2=\rho\eps.
\ee

\begin{figure}
\begin{center}
 {\bf a}\includegraphics[height=1.8in,width=2.0in]{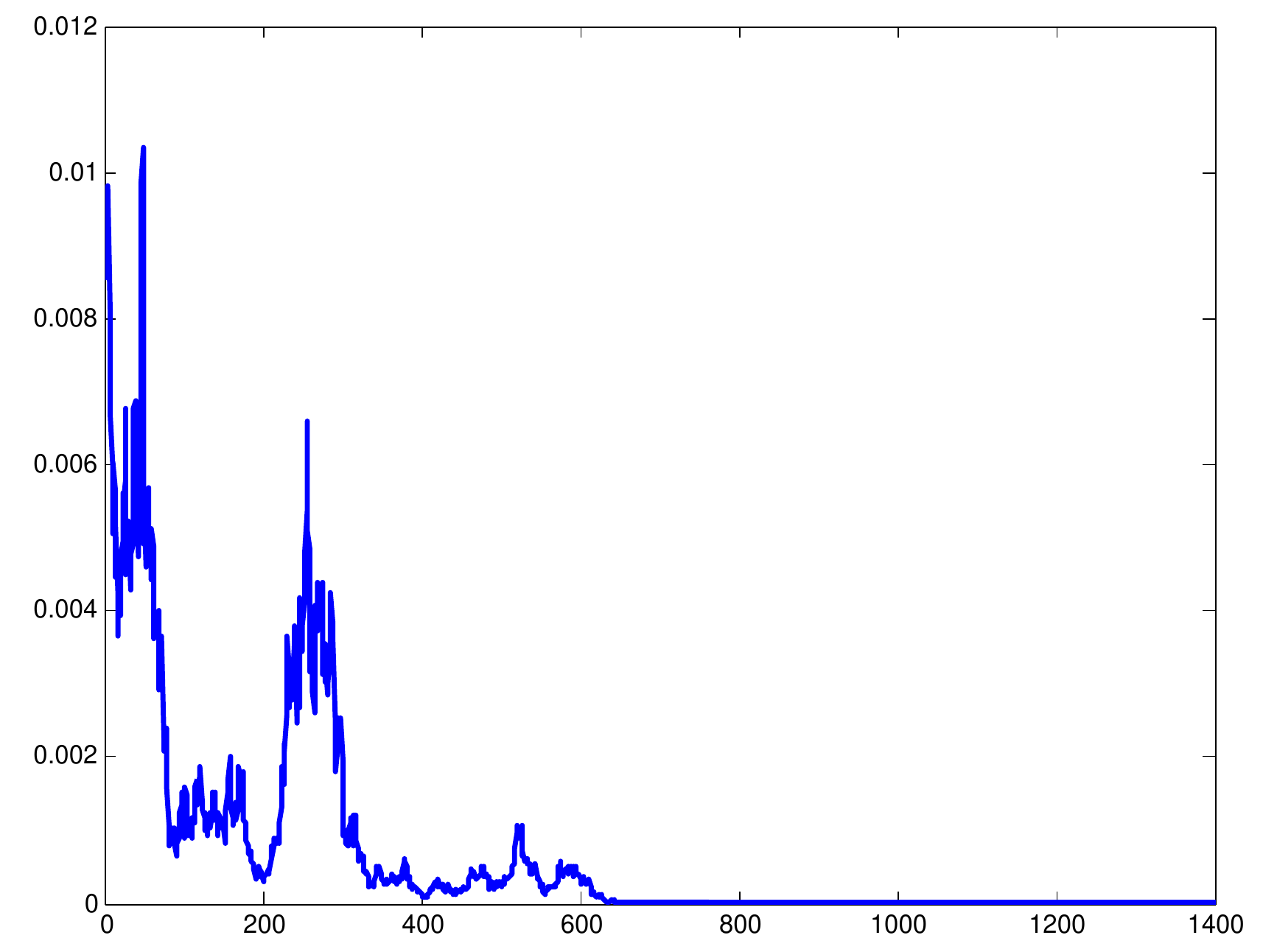}\quad
 {\bf b}\includegraphics[height=1.8in,width=2.0in]{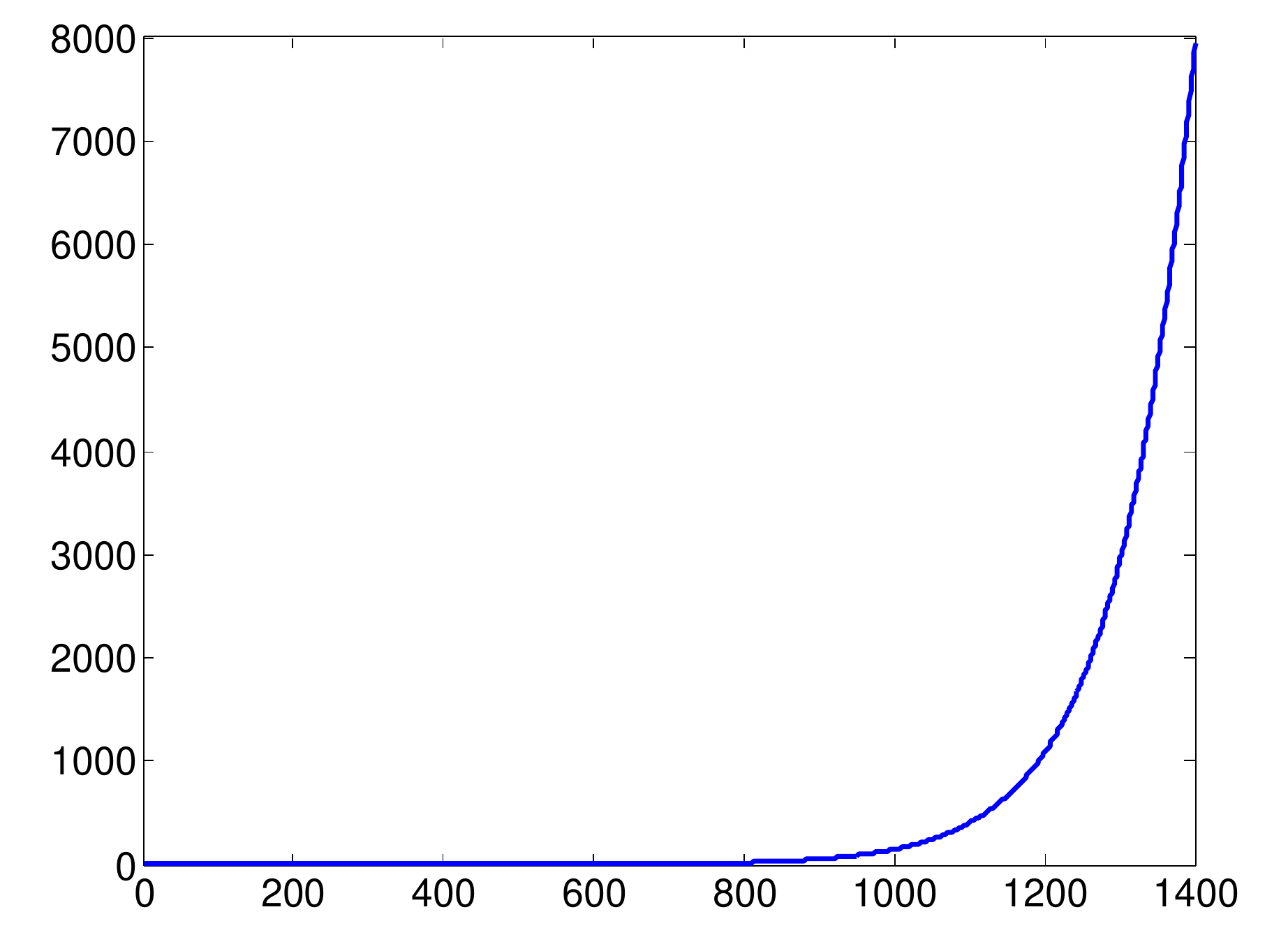}
\end{center}
\caption{ \textbf{a}) Time series $|x_n|$ generated by the  stochastic  two-dimensional 
system defined in 
Example~\ref{ex.2d-1}. The values of parameters are $\eps=0.01,$ $\rho=5$.
\textbf{b}) The time series generated by the underlying deterministic system
($\sigma=0$) is included for comparison.
}
\lbl{f.3}
\end{figure}

\begin{ex}\lbl{ex.2d-1} Consider (\ref{planar}) with matrix $B$ defined in (\ref{defB}) and
$$
A=\begin{pmatrix} 1+\eps & 0 \\ 0 & 0.5 \end{pmatrix}, \quad 0<\eps\ll 1.
$$
Figure~\ref{f.3}\textbf{a} shows a typical trajectory of the randomly perturbed system.
The noise keeps the trajectory from diverging from the origin which takes place in the
deterministic system (Figure~\ref{f.3}\textbf{b}).
\end{ex}

\begin{ex}\lbl{ex.2d-2} In this example, we consider a nonnormal matrix with multiple 
eigenvalues
$$
A=\begin{pmatrix} 1+\eps & 0.1 \\ 0 & 1+\eps \end{pmatrix}, \quad 0<\eps\ll 1.
$$
Figure~\ref{f.4} shows the results of the stabilization by noise for this case. The experiments
with the noise intensity in plots \textbf{a} and \textbf{b} show that stronger (albeit small)
noise results in a more robust stabilization.
\end{ex}

\vskip 0.2cm
\noindent
{\bf Acknowledgements.}
This work was supported in part by  a grant from Simons Foundation (grant \#208766 to
PH) and by the NSF (grant DMS 1412066 to GM). GM has also benefitted
from participating in a SQuaRe group `Stochastic stabilisation of
limit-cycle dynamics in ecology and neuroscience' sponsored by the
American Institute of Mathematics.

\begin{figure}
\begin{center}
 {\bf a}\includegraphics[height=1.8in,width=2.0in]{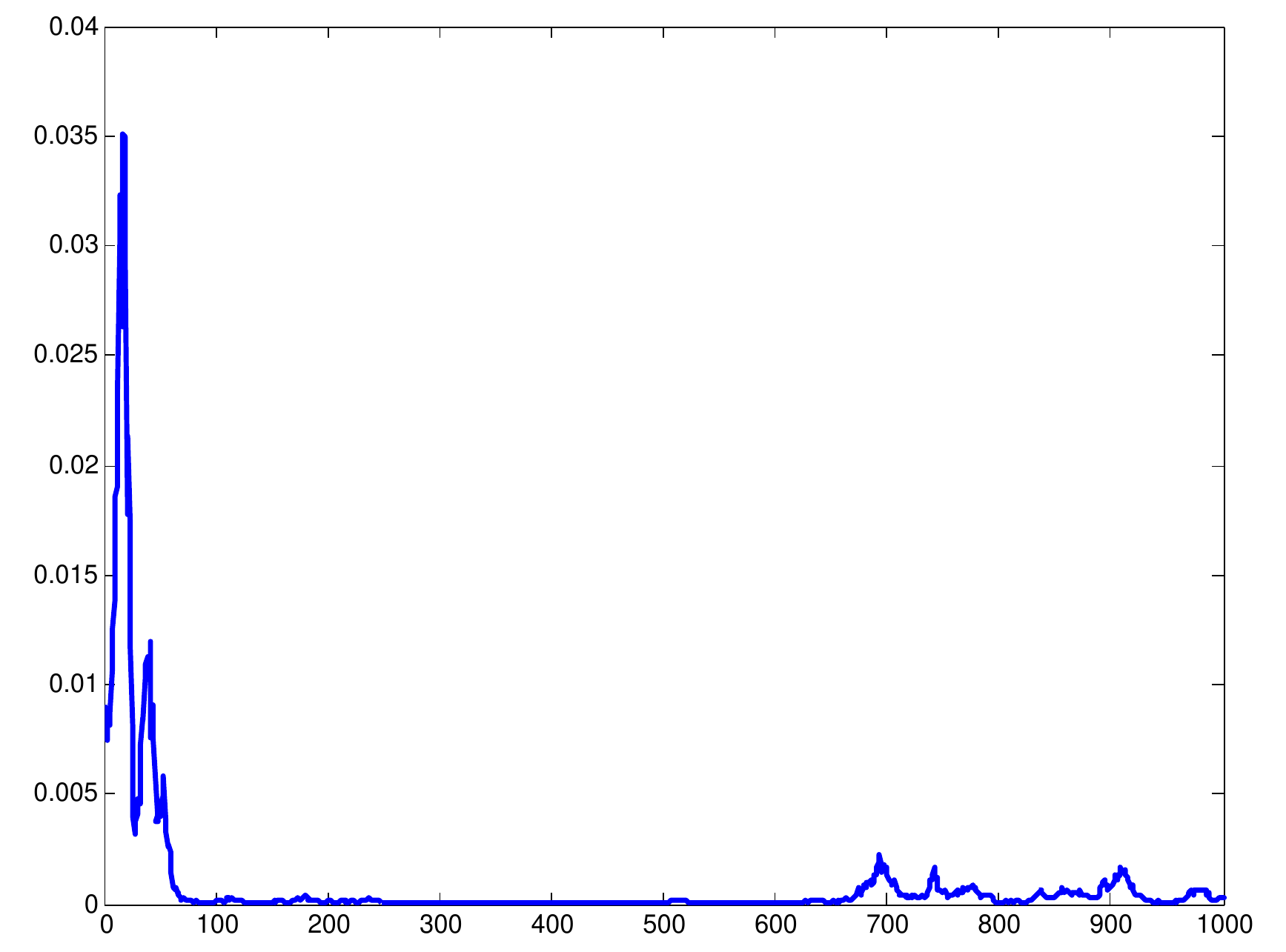}
 {\bf a}\includegraphics[height=1.8in,width=2.0in]{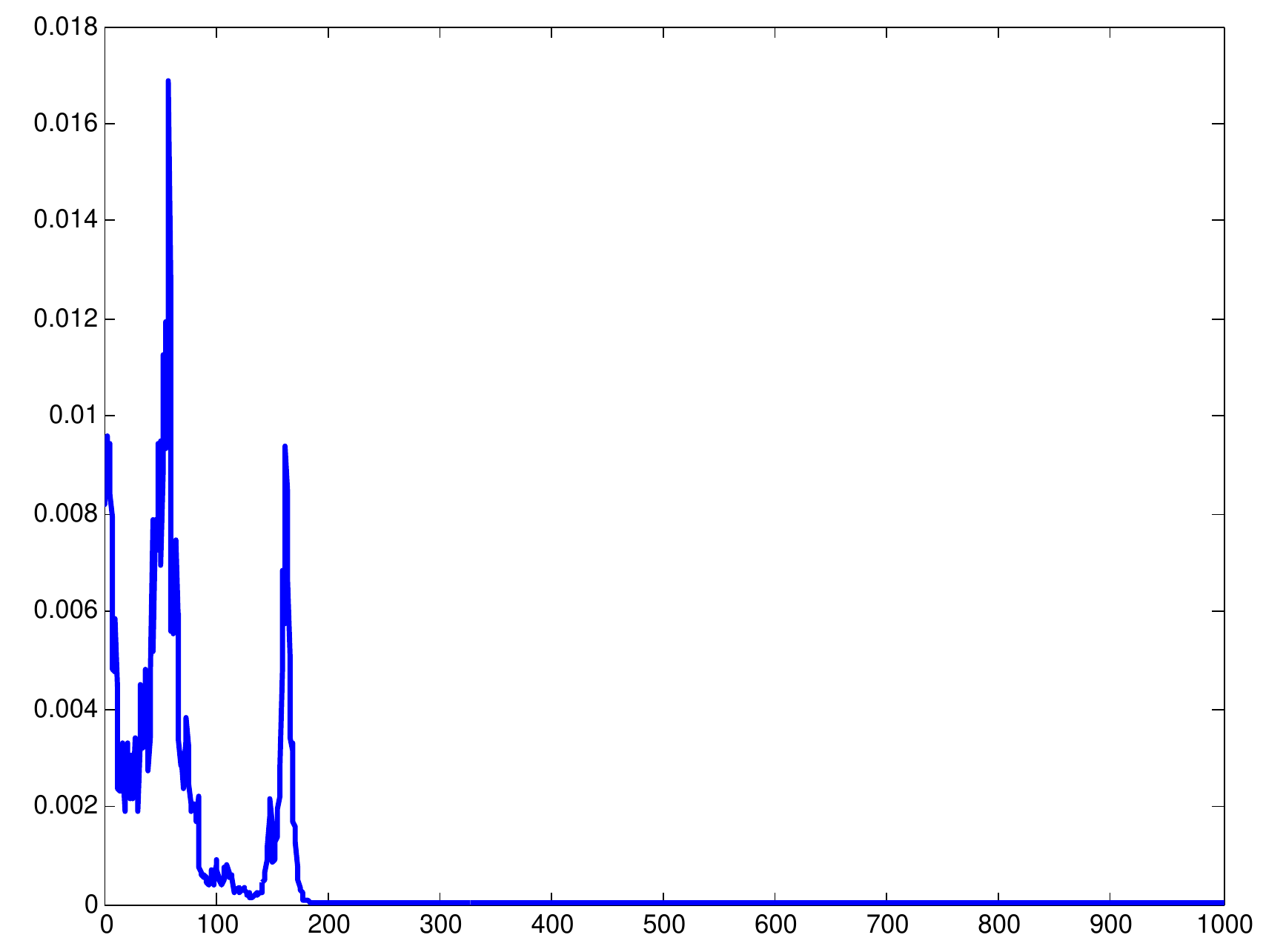}
 {\bf b}\includegraphics[height=1.8in,width=2.0in]{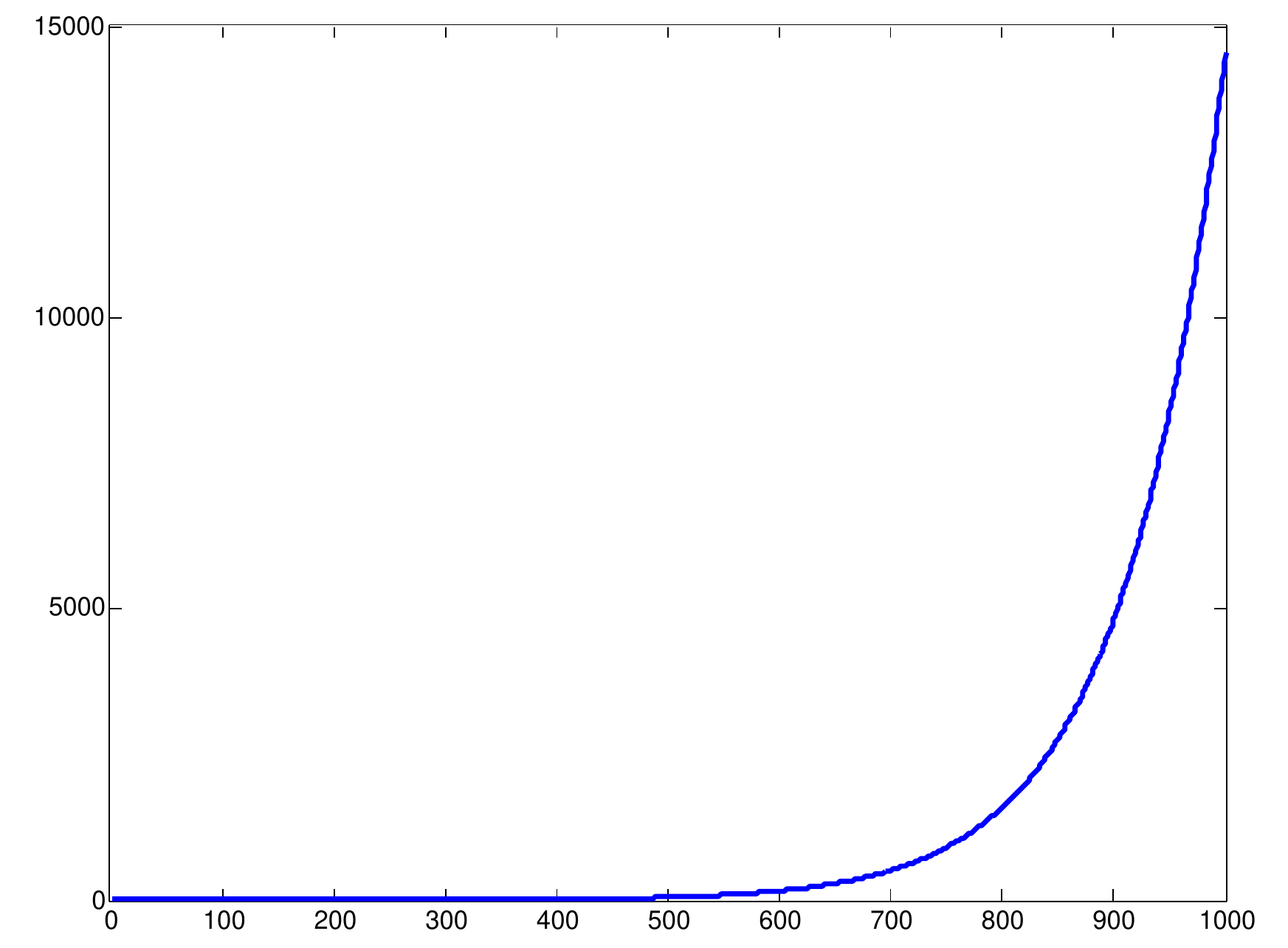}
\end{center}
\caption{ \textbf{a}) Time series  $|x_n|$ generated by the  stochastic  two-dimensional system defined in 
Example~\ref{ex.2d-2}. The values of parameters are $\eps=0.01,$ $\rho=5$.  \textbf{b})  
The same as in \textbf{a} but with $\rho=10$.
\textbf{c}) The time series generated by the underlying deterministic system
($\sigma=0$) is included for comparison.
}
\lbl{f.4}
\end{figure}

 \bibliographystyle{amsplain}

\begin{thebibliography}{10}

\bibitem{AVR86}
V.M.~Afra{\u\i}movich, N.N.~Verichev, and M.~I. Rabinovich, \emph{Stochastic
  synchronization of oscillations in dissipative systems}, Izv. Vyssh. Uchebn.
  Zaved. Radiofiz. \textbf{29} (1986), no.~9, 1050--1060. \MR{877439
  (88g:58110)}

\bibitem{ABR08}
J.~Appleby, G.~Berkolaiko, and A.~Rodkina, \emph{On local stability for a
  nonlinear difference equation with a non-hyperbolic equilibrium and fading
  stochastic perturbations}, J. Difference Equ. Appl. \textbf{14} (2008),
  no.~9, 923--951. \MR{2439782 (2009g:39007)}

\bibitem{AKMR10}
J.~Appleby, C.~Kelly, X.~Mao, and A.~Rodkina, \emph{On
  the local dynamics of polynomial difference equations with fading stochastic
  perturbations}, Dyn. Contin. Discrete Impuls. Syst. Ser. A Math. Anal.
  \textbf{17} (2010), no.~3, 401--430. \MR{2656407 (2011d:39028)}

\bibitem{ABR09}
J.~Appleby, G.~Berkolaiko, and A.~Rodkina,
  \emph{Non-exponential stability and decay rates in nonlinear stochastic
  difference equations with unbounded noise}, Stochastics \textbf{81} (2009),
  no.~2, 99--127. \MR{2571683 (2010j:39028)}

\bibitem{AppMao05}
J.~Appleby and X.~Mao, \emph{Stochastic stabilisation of
  functional differential equations}, Systems Control Lett. \textbf{54} (2005),
  no.~11, 1069--1081. \MR{2170288 (2006d:34156)}

\bibitem{AMR06}
J.~Appleby, X.~Mao, and A.~Rodkina, \emph{On stochastic
  stabilization of difference equations}, Discrete Contin. Dyn. Syst.
  \textbf{15} (2006), no.~3, 843--857. \MR{2220752 (2007b:39031)}

\bibitem{Arn90}
L.~Arnold, \emph{Stabilization by noise revisited}, Z. Angew. Math. Mech.
  \textbf{70} (1990), no.~7, 235--246. \MR{1066866 (91j:93119)}

\bibitem{BerGen06}
N.~Berglund and B.~Gentz, \emph{Noise-induced phenomena in slow-fast
  dynamical systems}, Probability and its Applications (New York),
  Springer-Verlag London, Ltd., London, 2006, A sample-paths approach.
  \MR{2197663 (2007b:37115)}

\bibitem{BR06}
G.~Berkolaiko and A.~Rodkina, \emph{Almost sure convergence of
  solutions to nonhomogeneous stochastic difference equation}, J. Difference
  Equ. Appl. \textbf{12} (2006), no.~6, 535--553. \MR{2240374 (2007b:39002)}

\bibitem{bil}
P.~Billingsley, \emph{Probability and measure}, 3rd ed., Wiley, 1995.

\bibitem{BR12}
E.~Braverman and A.~Rodkina, \emph{On difference equations with asymptotically
  stable 2-cycles perturbed by a decaying noise}, Comput. Math. Appl.
  \textbf{64} (2012), no.~7, 2224--2232. \MR{2966858}

\bibitem{BucKel10}
E.~Buckwar and C.~Kelly, \emph{Towards a systematic linear
  stability analysis of numerical methods for systems of stochastic
  differential equations}, SIAM J. Numer. Anal. \textbf{48} (2010), no.~1,
  298--321. \MR{2608371 (2011b:60271)}

\bibitem{DVM05}
R.E.L.~DeVille, E.~Vanden-Eijnden, and C.B.~Muratov, \emph{Two
  distinct mechanisms of coherence in randomly perturbed dynamical systems},
  Phys. Rev. E (3) \textbf{72} (2005), no.~3, 031105, 10. \MR{2179903
  (2006f:37074)}

\bibitem{DRR06}
B.~Doiron, J.~Rinzel, and A.~Reyes, \emph{Stochastic synchronization in
  finite size spiking networks}, Phys. Rev. E (3) \textbf{74} (2006), no.~3,
  030903, 4. \MR{2282117 (2007k:92017)}

\bibitem{Fre01}
M.~Freidlin, \emph{On stochastic perturbations of dynamical systems with fast
  and slow components}, Stoch. Dyn. \textbf{1} (2001), no.~2, 261--281.
  \MR{1840196 (2003a:60032)}

\bibitem{FurKes60}
H.~Furstenberg and H.~Kesten, \emph{Products of random matrices}, Ann. Math.
  Statist. \textbf{31} (1960), 457--469.

\bibitem{GP05}
D.S.~Goldobin and A.~Pikovsky, \emph{Synchronization and
  desynchronization of self-sustained oscillators by common noise}, Phys. Rev.
  E (3) \textbf{71} (2005), no.~4, 045201, 4. \MR{2139983 (2005m:82085)}

\bibitem{Hig00}
D.J.~Higham, \emph{Mean-square and asymptotic stability of the stochastic
  theta method}, SIAM J. Numer. Anal. \textbf{38} (2000), no.~3, 753--769
  (electronic). \MR{1781202}

\bibitem{HMY07}
D.J.~Higham, X.~Mao, and C.~Yuan, \emph{Almost sure and moment
  exponential stability in the numerical simulation of stochastic differential
  equations}, SIAM J. Numer. Anal. \textbf{45} (2007), no.~2, 592--609
  (electronic). \MR{2300289 (2008c:60064)}

\bibitem{HM09}
P.~Hitczenko and G.S.~Medvedev, \emph{Bursting oscillations induced by
  small noise}, SIAM J. Appl. Math. \textbf{69} (2009), no.~5, 1359--1392.
  \MR{2487064 (2010f:60169)}

\bibitem{HM13}
\bysame, \emph{The {P}oincar\'e map of randomly perturbed periodic motion}, J.
  Nonlinear Sci. \textbf{23} (2013), no.~5, 835--861. \MR{3101836}

\bibitem{Horn-Johnson}
R.A.~Horn and C.R.~Johnson, \emph{Matrix analysis}, second ed.,
  Cambridge University Press, Cambridge, 2013. \MR{2978290}

\bibitem{Kha-StochStability}
R.~Khasminskii, \emph{Stochastic stability of differential equations},
  second ed., Stochastic Modelling and Applied Probability, vol.~66, Springer,
  Heidelberg, 2012, With contributions by G. N. Milstein and M. B. Nevelson.
  \MR{2894052}

\bibitem{KocPal10}
H.~Ko{\c{c}}ak and K.J.~Palmer, \emph{Lyapunov exponents and
  sensitivity dependence}, J. Dynam. Differential Equations \textbf{22} (2010),
  no.~3, 381--398. \MR{2719912 (2012f:37075)}

\bibitem{LaiLor10}
C.~Laing and G.J.~Lord (eds.), \emph{Stochastic methods in
  neuroscience}, Oxford University Press, Oxford, 2010. \MR{2640514
  (2010m:60006)}

\bibitem{Long10}
A.~Longtin, \emph{Neural coherence and stochastic resonance}, Stochastic
  methods in neuroscience, Oxford Univ. Press, Oxford, 2010, pp.~94--123.
  \MR{2642697}

\bibitem{Mao94}
X.~Mao, \emph{Stochastic stabilization and destabilization}, Systems
  Control Lett. \textbf{23} (1994), no.~4, 279--290. \MR{1298174 (95h:93089)}

\bibitem{MP08}
M.~Porfiri and R.~Pigliacampo, \emph{Master-slave global stochastic
  synchronization of chaotic oscillators}, SIAM J. Appl. Dyn. Syst. \textbf{7}
  (2008), no.~3, 825--842. \MR{2443024 (2009h:93117)}

\bibitem{SaiMit96}
Y.~Saito and T.~Mitsui, \emph{Stability analysis of numerical
  schemes for stochastic differential equations}, SIAM J. Numer. Anal.
  \textbf{33} (1996), no.~6, 2254--2267. \MR{1427462 (98c:65138)}

\bibitem{SaiMit02}
\bysame, \emph{Mean-square stability of numerical schemes for stochastic
  differential systems}, Vietnam J. Math. \textbf{30} (2002), no.~suppl.,
  551--560. \MR{1964242 (2003m:65014)}
\end{thebibliography}

\def\cprime{$'$}
\providecommand{\bysame}{\leavevmode\hbox to3em{\hrulefill}\thinspace}
\providecommand{\MR}{\relax\ifhmode\unskip\space\fi MR }
\providecommand{\MRhref}[2]{%
  \href{http://www.ams.org/mathscinet-getitem?mr=#1}{#2}
}
\providecommand{\href}[2]{#2}

\end{document}